\documentclass{article}

\usepackage{amsfonts}
\usepackage{amssymb}

\date{}
\title{\Large \bf Bornologies, selection principles and\\
 function spaces}

\author{A. Caserta\footnote{Supported by GNSAGA}, \, G. Di Maio\footnote{Supported
by GNSAGA}, \, Lj.D.R. Ko\v cinac\footnote{Supported by MNTR RS,
GNSAGA and SUN} }

{}
\newtheorem{theorem}{{\bf Theorem}}[section]

\newtheorem{corollary}[theorem]{{\bf Corollary}}

\newcommand{\proof}{{\bf Proof.  }}
\newcommand{\eproof}{{\blacktriangle}}

\newcommand{\sone}{{\sf S}_1}
\newcommand{\sfin}{{\sf S}_{fin}}

\newcommand{\naturals}{{\mathbb N}}
\newcommand{\reals}{{\mathbb R}}
\newcommand{\zero}{\underline{\sf 0}}

\begin{document}

\maketitle

\begin{abstract}
We study some closure-type properties of function spaces endowed
with the new topology of strong uniform convergence on a bornology
introduced by Beer and Levy in 2009. The study of these function
spaces was initiated in \cite{beerlevi} and \cite{agataholapeppe}.
The properties we study are related to selection principles.
\end{abstract}

\medskip
\begin{flushleft}
\vspace{.3cm} {\sf Keywords}: Bornology, topology of strong uniform
convergence, $\mathfrak B^s$-cover, selection principles, function
spaces.
\end{flushleft}

\section{Introduction}

Our notation and terminology are standard as in \cite{engelking}.
All spaces are assumed to be metric.

If $(X,d)$ is a metric space, $x\in X$, $A\subset X$ and
$\varepsilon >0$ a real number, we write
\begin{center}
$S(x,\varepsilon)= \{y\in X:d(x,y) <\varepsilon\}$,\\
$A^{\varepsilon} := \bigcup_{a\in A}S(a,\varepsilon)$,
\end{center}
to denote the open $\varepsilon$-ball with center $x$ and the
$\varepsilon$-enlargement of $A$.

Given spaces $X$ and  $Y$ we denote by $Y^X$ (resp. ${\sf C}(X,Y)$) the set
of all functions (resp. all continuous functions) from $X$ into $Y$. When $Y =
\reals$ we write ${\sf C}(X)$ instead of ${\sf C}(X,\reals)$. ${\sf C}_p(X)$ and
${\sf C}_k(X)$ denote the set ${\sf C}(X)$ endowed with the pointwise topology $\tau_p$ and
the compact-open topology $\tau_k$, respectively.

Recall that a \emph{bornology} on a metric space $(X,d)$ is a
family $\mathfrak B$ of nonempty subsets of $X$ which is closed
under finite unions, hereditary (i.e. closed under taking nonempty
subsets) and forms a cover of $X$ \cite{hu}, \cite{hogbe}.
Throughout the paper we suppose that $X$ does not belong to a
bornology $\mathfrak B$ on $X$. A \emph{base} for a bornology
$\mathfrak B$ on $(X,d)$ is a subfamily $\mathfrak B_0$ of
$\mathfrak B$ which is cofinal in $\mathfrak B$ with respect to
the inclusion, i.e. for each $B\in\mathfrak B$ there is
$B_0\in\mathfrak B_0$ such that $B\subset B_0$. A base is called
\emph{closed} (\emph{compact})  if all its members are closed
(compact) subsets of $X$.

For example, if $(X,d)$ is a metric space, then the family $\mathfrak
F$ of all nonempty finite subsets of $X$ is a bornology on $X$;
it is the smallest bornology on $X$ and has a closed (in fact a compact)
base. The largest bornology on $X$ is the set of all nonempty subsets of $X$.
Other important bornologies are: the collection of all nonempty (i)
relatively compact subsets (i.e. subsets with compact closure), denoted by
$\mathfrak K$,
(ii) $d$-bounded subsets, (iii) totally $d$-bounded subsets.

The following simple facts will be used in the sequel.

\noindent {\bf Fact 1.} If a bornology $\mathfrak B$ has a closed
base, then $B \in \mathfrak B$ implies $\overline{B} \in \mathfrak
B$.

\noindent {\bf Fact 2.} For every $B$ in a bornology $\mathfrak B$ and every
$\delta > 0$ it holds $\overline{B^{\delta}} \subset B^{2\delta}$.

In \cite{beerlevi} the notion of strong uniform continuity was
introduced: a (not necessarily continuous) mapping $f:X \to Y$ from
a metric space $(X,d)$ to a metric space $(Y,\rho)$ is
\emph{strongly uniformly continuous} on a subset $B$ of $X$ if for
every $\varepsilon >0$ there is $\delta > 0$ such that $d(x_1,x_2) <
\delta$ and $\{x_1,x_2\} \cap B \neq \emptyset$ imply
$\rho(f(x_1),f(x_2)) < \varepsilon$. If $\mathfrak B$ is a bornology
on $X$, then $f:X\to Y$ is called \emph{strongly uniformly
continuous on $\mathfrak B$} if $f$ is strongly uniformly continuous
on $B$ for each $B\in \mathfrak B$.

In the same paper Beer and Levi defined a new topology on the set
$Y^X$ of all functions from $X$ into $Y$, named the topology of
strong uniform convergence. They initiated the study of function
spaces $Y^X$ and ${\sf C}(X,Y)$ with this new topology and
characterized metrizability and first countability. Further
analysis of these function spaces was done in
\cite{agataholapeppe} for several topological properties (for
example, (sub)metrizability, first countability, complete
metrizabilty, separability, countable tightness, the
Fr\'echet-Urysohn property). We continue the investigation in this
area considering properties related to selection principles. Also,
we indicate how the idea of strong uniform convergence can be
applied to general selection principles theory. For more
information about selection principles see the survey papers
\cite{iransurv}, \cite{marionsurv}, \cite{boazsurv} and references
therein.

If $\mathfrak U$ and  $\mathfrak V$ are collections of subsets of a
space $X$. Then:

(1) $\sone({\mathfrak U},{\mathfrak
V})$ denotes the selection hypothesis that for each sequence
$(U_n:n\in\naturals)$ of elements of ${\mathfrak U}$ there is a
sequence $(u_n:n\in\naturals)$ such that for each
$n\in\naturals$, $u_n\in U_n$ and $\{u_n:n\in\naturals\}$ is in
$\mathfrak V$;

(2) $\sfin({\mathfrak U},{\mathfrak V})$ is the selection hypothesis
that for each sequence $(U_n:n\in\naturals)$ of elements of ${\mathfrak U}$
there is a sequence $(V_n:n\in\naturals)$ such that $V_n$ is a
finite subset of $U_n$ for each $n\in\naturals$ and
$\bigcup_{n\in\mathbb N}V_n \in {\mathfrak V}$.


\section{Function spaces}


We begin with a definition from \cite{beerlevi}.

For given metric spaces $(X,d)$ and $(Y,\rho)$ and a bornology
$\mathfrak B$ {\bf with closed base} on $X$ we denote by $\tau_{\mathfrak
B}^s$ the \emph{topology of strong uniform convergence on $\mathfrak B$}
determined by a uniformity on $Y^X$ having as a base the sets of the form
\[
[B,\varepsilon]^s := \{(f,g): \exists \delta >0 \,\forall  x\in
B^{\delta}, \rho(f(x),g(x)) < \varepsilon\} \, \, \, (B\in
\mathfrak B, \varepsilon >0).
\]
Therefore, for a function $f\in ({\sf C}(X,Y),\tau_{\mathfrak B}^s)$
the standard local base of $f$ is the collection of sets
\[
[B,\varepsilon]^s(f) =\{g\in ({\sf C}(X,Y),\tau_{\mathfrak
B}^s):\exists \delta > 0, \, \rho(g(x),f(x))<\varepsilon, \,
\forall x\in B^{\delta}\} \, \, \, (B\in \mathfrak B, \varepsilon
>0).
\]

For each bornology $\mathfrak B$ with closed base on $X$ the
topology $\tau_{\mathfrak B}^s$ on $Y^X$ is finer than the
classical topology $\tau_{\mathfrak B}$ of uniform convergence on
$\mathfrak B$, and if $\mathfrak B$ has a compact base, then
$\tau_{\mathfrak B}^s = \tau_{\mathfrak B}\le \tau_k$ on ${\sf
C}(X,Y)$. In particular, $\tau_p \le \tau_{\mathfrak F}^s \le
\tau_{\mathfrak B}^s \le \tau_{\mathfrak K}^s=\tau_k$ on ${\sf
C}(X)$ (see \cite[Corollary 6.6]{beerlevi}, \cite[Remark
3.4]{agataholapeppe}).

\medskip
The following notion was introduced in \cite{agataholapeppe}. An
open cover $\mathcal U$ of a metric space $(X,d)$ with a bornology
$\mathfrak B$ is said to be a \emph{strong $\mathfrak B$-cover} of
$X$ (or shortly a \emph{$\mathfrak B^s$-cover} of $X$) if $X\notin
\mathcal U$ and for each $B\in \mathfrak B$ there exist $U\in
\mathcal U$ and $\delta > 0$ such that $B^{\delta} \subset U$.

In this paper the collection of all strong $\mathfrak B$-covers of a
space is denoted by $\mathcal O_{\mathfrak B^s}$. We also suppose that all
spaces we consider are \emph{$\mathfrak B^s$-Lindel\"of}, i.e. each $\mathfrak B^s$-cover
contains a countable $\mathfrak B^s$-subcover.

Let us define also the following. A countable open cover $\mathcal U =
\{U_n:n\in \naturals\}$ of $X$ is said to be a \emph{$\gamma_s$-cover} (\emph{$\gamma_{\mathfrak
B^s}$-cover}, called a \emph{$\mathfrak B^s$-sequence} in \cite{agataholapeppe})
if it is infinite and for $x\in X$ (each $B\in \mathfrak B$) there are $n_0\in\naturals$
and a sequence $(\delta_n:n\ge n_0)$ of positive real numbers such that $S(x,\delta_n)\subset U_n$ ($B^{\delta_n}
\subset U_n)$ for all $n\ge n_0$.

For a given bornology
$\mathfrak B$ in a space $X$ we denote by $\Gamma_s$ and
$\Gamma_{\mathfrak B^s}$ the collection of all countable $\gamma_s$-covers
and the collection of all $\gamma_{\mathfrak B^s}$-covers of $X$.

The symbol $\zero$ denotes the constantly zero function in $({\sf C}(X),
\tau_{\mathfrak B}^s)$. The space $({\sf C}(X),\tau_{\mathfrak B}^s)$ is homogeneous so that
it suffices to look at the point $\zero$ when studying local
properties of this space.

In \cite{agataholapeppe} the following theorem was proved.

\begin{theorem} \label{tightbess} Let $(X,d)$ be a metric space and $\mathfrak B$ a
bornology on $X$ with closed base. The following are equivalent:
\begin{itemize}
\item[$(1)$] $({\sf C}(X),\tau_{\mathfrak B}^s)$ has countable
tightness;
\item[$(2)$] $X$ is a $\mathfrak B^s$-Lindel\"of space.
\end{itemize}
\end{theorem}

Our aim is to prove similar results for countable fan tightness
and countable strong fan tightness.

For a space $X$ and a point $x\in X$ the symbol $\Omega_x$ denotes
the set $\{A\subset X\setminus\{x\}: x\in \overline{A}\}$, and $\Sigma_x$ is the set of sequences
converging to $x$.

A space $X$ has \emph{countable fan tightness} \cite{arhbook} if for each $x\in X$
we have that $\sfin(\Omega_x,\Omega_x)$ holds. $X$ has {\it
countable strong fan tightness} \cite{sakai} if for each $x\in X$ the selection
principle $\sone(\Omega_x,\Omega_x)$ holds.

\begin{theorem} \label{strongfan} Let $(X,d)$ be a metric space and $\mathfrak B$ a
bornology on $X$ with closed base. The following are equivalent:
\begin{itemize}
\item[$(1)$] $({\sf C}(X),\tau_{\mathfrak B}^s)$ has countable
strong fan tightness;
\item[$(2)$] $X$ satisfies $\sone(\mathcal O_{\mathfrak
B^s},\mathcal O_{\mathfrak B^s})$.
\end{itemize}
\end{theorem}
$\proof$ $(1) \Rightarrow (2)$: Let $(\mathcal U_n:n\in\naturals)$
be a sequence of open strong $\mathfrak B$-covers of $X$. For every
$n\in \naturals$ and every $B\in \mathfrak B$ there exists $U_B \in
\mathcal U_n$ and $\delta >0$ such that $B^{2\delta} \subset U_B$.

For a fixed $n\in\naturals$ and $B\in\mathfrak B$ let
\[
\mathcal{U}_{n,B}:=\{U\in{\mathcal U} _n: B^{2\delta}\subset U\}.
\]
For each $U\in\mathcal{U}_{n,B}$, by Fact 2, there is a continuous
function $f_{B,U}$ from $X$ into $[0,1]$ such that
$f_{B,U}(B^{\delta})= \{0\}$ and $f_{B,U}(X\setminus U)=\{1\}$. Let
for each $n$,
\[
A_n=\{f_{B,U}:B\in  \mathfrak B, U\in \mathcal{U}_{n,B}\}.
\]
It is easily seen that $\zero$ belongs to the closure of $A_n$ in
$\tau_{\mathfrak B}^s$ for each $n\in\naturals$. By (1) there is a
sequence $(f_{B_n,U_n}:n\in\naturals)$ such that for each $n$,
$f_{B_n,U_n} \in A_n$ and $\zero\in\overline{\{f_{B_n,U_n}:n\in
\naturals\}}$. We claim that the sequence $(U_n:n\in\naturals)$
witnesses that $X$ has property $\sone(\mathcal O_{\mathfrak
B^s},\mathcal O_{\mathfrak B^s})$.

Let $B\in \mathfrak B$. Since $\zero\in\overline{\{f_ {B_n,U_n}:n\in
\naturals\}}$ it follows that there is $m\in\naturals$ such that
$[B,1]^s(\zero)$ contains the function $f_{B_m,U_m}$. Therefore,
there is $\delta > 0$ such that for each $x\in B^{\delta}$ it holds
$f_{B_m,U_m}(x) <1$, which means $B^{\delta} \subset U_m$.

\medskip
$(2)\Rightarrow (1)$: \ Let $(A_n:n\in\naturals)$ be a
sequence of subsets of $({\sf C}(X),\tau_{\mathfrak B}^s)$ whose
closures contain $\zero$.

For every $B\in \mathfrak B$ and every $m\in \naturals$ the
neighborhood $[B,1/m]^s(\zero)$ of $\zero$ intersects each $A_n$.
It follows that for each $n\in \naturals$ there exists a function
$f_{B,n,m}\in A_n$ satisfying: there is $\delta > 0$ with
$|f_{B,n,m}(x)|< 1/m$ for each $x\in B^{\delta}$. For each $n$ set
\[
{\mathcal U}_{n,m}=\{f^{\gets}(-1/m,1/m):m\in \naturals, f\in A_n\}.
\]
(We can view the indices $m,n$ in $\mathcal U_{n,m}$ as
$\varphi(m,n)$ for some bijection  $\varphi:\naturals^2 \to
\naturals$.) We claim that for each $n,m\in \naturals$, ${\mathcal
U}_{n,m}$ is a $\mathfrak B^s$-cover of $X$. Indeed, if $B\in
\mathfrak B$, then there is $f_{B,n,m}\in [B,1/m]^s(\zero) \cap
A_n$. Hence there is $\delta
> 0$ such that $|f_{B,n,m}(x)| < 1/m$ for each $x\in B^{\delta}$. This means
$B^{\delta} \subset f_{B,n,m}^{\gets}(-1/m,1/m) \in \mathcal
U_{n,m}$.

Put $M = \{m\in\naturals: X \in \mathcal U_{n,m} \mbox{ for some }
n\in \naturals\}$.

\smallskip
Case 1. $M$ is infinite.

\smallskip
There are $m_1 < m_2 < \cdots$ in $M$ and (the corresponding) $n_1,
n_2, \cdots$ in $\naturals$ such that $f_{B_i,n_i,m_i}^{\gets}
(-1/m_i,1/m_i) = X$ for all $i\in \naturals$ and some $B_i\in
\mathfrak B$. Let $[B,\varepsilon]^s(\zero)$ be a $\tau_{\mathfrak
B}^s$-neighbourhood of $\zero$. Pick $m_k$ such that $1/m_k <
\varepsilon$. For every $m_i > m_k$ we have $f_{B_i,n_i,m_i}(x) \in
(-1/m_i,1/m_i)$ for each $x \in X$ and so $f_{B_i,n_i,m_i} \in
[B,1/m_i]^s(\zero) \subset [B,\varepsilon]^s(\zero)$. This means
that the sequence $(f_{B_i,n_i,m_i}:i\in \naturals)$ converges to
$\zero$, hence $\zero \in \overline{\{f_{B_i,n_i,m_i}:i\in
\naturals\}}$.

\smallskip
Case 2. $M$ is finite.

\smallskip
There is $m_0\in \naturals$ such that for each $m\ge m_0$ and each
$n\in \naturals$, the set $\mathcal U_{n,m}$ is a $\mathfrak
B^s$-cover of $X$. One may suppose $m_0=1$, and since the set $
\{n\in\naturals: X \in \mathcal U_{n,n}\}$ is also finite we can
work with this set instead of $M$, and assume that for each $n\in
\naturals$, $\mathcal U_{n,n}$ is a $\mathfrak B^s$-cover of $X$.
By (2) choose for each $n\in \naturals$ a set $U_{n,n}\in \mathcal
U_{n,n}$ so that $\{U_{n,n}:n\in \naturals\}$ is a $\mathfrak
B^s$-cover of $X$. Consider the corresponding functions
$f_{B_n,n,n}$, $n\in \naturals$. We claim that $\zero \in
\overline{\{f_{B_n,n,n}:n\in \naturals\}}$. Let
$[B,\varepsilon]^s(\zero)$ be a neighbourhood of $\zero$. There
are $j\in \naturals$ and $\delta > 0$ such that $B^\delta \subset
U_{j,j}$. But the set $K$ of all such $j$ is infinite because
$\{U_{n,n}:n\in \naturals\}$ is a $\mathfrak B^s$-cover of $X$.
Take $k \in K$ so that $1/k < \varepsilon$. Then $B^\delta \subset
f_{B_{k},k,k}^\gets(-1/k,1/k) \subset
f_{B_{k},k,k}^\gets(-\varepsilon,\varepsilon)$, i.e.
$f_{B_{k},k,k} \in [B,\varepsilon]^s(\zero)$. $\eproof$

\medskip
Let $\mathcal K$ denote the family of $k$-covers of a space $X$. (An
open cover is a $k$-cover if each compact
set $K\subset X$ is contained in a member of the cover.)

\begin{corollary}\label{k-strongfan} ({\rm (\cite{agt})} The space ${\sf C}_k(X)$
has countable strong fan tightness if and only if $X$ has property
$\sone (\mathcal K, \mathcal K)$.
\end{corollary}

\medskip
The following theorem can be proved similarly.

\begin{theorem} \label{fan} Let $(X,d)$ be a metric space and $\mathfrak B$ a
bornology with closed base on $X$. The following are equivalent:
\begin{itemize}
\item[$(1)$] $({\sf C}(X),\tau_{\mathfrak B}^s)$ has countable
fan tightness;
\item[$(2)$] $X$ satisfies $\sfin(\mathcal O_{\mathfrak
B^s},\mathcal O_{\mathfrak B^s})$.
\end{itemize}
\end{theorem}

\begin{corollary}\label{k-fan} ({\rm (\cite{chinese}, \cite{agt})} The space ${\sf C}_k(X)$
has countable fan tightness if and only if $X$ has property $\sfin
(\mathcal K, \mathcal K)$.
\end{corollary}

Recall that a space $X$ is said to be \emph{Fr\'echet-Urysohn} if
for each $A$ subset of $X$ and each $x\in \overline{A}$ there is a
sequence in $A$ converging to $x$. $X$ is \emph{strictly
Fr\'echet-Urysohn} if fulfills the selection property
$\sone(\Omega_x,\Sigma_x)$.

\begin{theorem} \label{SFU} Let $(X,d)$ be a metric space and $\mathfrak B$ be a
bornology on $X$. The following are equivalent:
\begin{itemize}
\item[$(1)$] $({\sf C}(X),\tau_{\mathfrak B}^s)$ is a strictly Fr\'echet-Urysohn space;
\item[$(2)$] $X$ satisfies $\sone(\mathcal O_{\mathfrak
B^s},\Gamma_{\mathfrak B^s})$.
\end{itemize}
\end{theorem}
$\proof$ $(1) \Rightarrow (2)$: Let $(\mathcal U_n:n\in\naturals)$
be a sequence of $\mathfrak B^s$-covers of $X$. For every $n\in
\naturals$ and every $B\in \mathfrak B$ there exist $U_{B,n} \in
\mathcal U_n$ and $\delta >0$ such that $B^{2\delta} \subset
U_{B,n}$. Set $\mathcal{U}_{n,B}:=\{U\in{\mathcal U}_{n}:
B^{2\delta}\subset U\}$. For each $U\in\mathcal{U}_{n,B}$ pick a
continuous function $f_{B,U}$ from $X$ into $[0,1]$ such that
$f_{B,U}(B^{\delta})= \{0\}$ and $f_{B,U}(X\setminus U)=\{1\}$.
Let for each $n$,
\[
A_n=\{f_{B,U}:B\in  \mathfrak B, U\in \mathcal{U}_{n,B}\}.
\]
Clearly the function $\zero$ belongs to the $\tau_{\mathfrak
B}^s$-closure of $A_n$ for each $n\in\naturals$, and since $({\sf
C}(X),\tau_{\mathfrak B}^s)$ is strictly Fr\'echet-Urysohn, there
are $f_{B_n,U_n} \in A_n$, $n\in \naturals$, such that the
sequence $(f_{B_n,U_n}:n\in\naturals)$ $\tau_{\mathfrak
B}^s$-converges to $\zero$. We prove that the set
$\{U_n:n\in\naturals\}$ is a $\gamma_{\mathfrak B^s}$-cover of
$X$. For the neighbourhood $[B,1]^s(\zero)$ of $\zero$ there is
$n_0\in \naturals$  such that for all $n > n_0$ $f_{B_n,U_n} \in
[B,1]^s(\zero)$. So, for each $n> n_0$ there is $\delta_n >0$ such
that $B^{\delta_n} \subset f_{B_n,U_n}^{\gets}(-1,1)$, hence
$B^{\delta_n} \subset U_n$.

\medskip
$(2) \Rightarrow (1)$: \ Let $(A_n:n\in\naturals)$ be a sequence of
subsets of $({\sf C}(X),\tau_{\mathfrak B}^s)$ such that $\zero \in
\bigcap_{n\in\naturals}\overline{A_n}$. For every $B\in \mathfrak B$
and every $m\in \naturals$ the neighborhood $[B,1/m]^s(\zero)$ of
$\zero$ intersects each $A_n$, and thus for each $n\in \naturals$
there is a function $f_{B,n,m}\in A_n$ such that there is $\delta
> 0$ with $|f_{B,n,m}(x)|< 1/m$ for each $x\in B^{\delta}$. For each
$n\in \naturals$ let
\[
{\mathcal U}_{n,m}=\{f^{\gets}(-1/m,1/m):m\in \naturals, f\in A_n\}.
\]
As in the proof of Theorem \ref{strongfan} we conclude that for each
$n,m\in \naturals$, ${\mathcal U}_{n,m}$ is a $\mathfrak B^s$-cover
of $X$. Apply now assumption (2) to the sequence $(\mathcal
U_{n,n}:n\in \naturals)$ and for each $n$ pick an element $U_{n,n}$
in $\mathcal U_{n,n}$ such that the set $\{U_{n,n}:n\in \naturals\}$
is a $\gamma_{\mathfrak B^s}$-cover of $X$. To each $U_{n,n}$
associate the corresponding function $f_{B_n,n,n} \in A_n$. We prove
that the sequence $(f_{B_n,n,n}:n\in\naturals)$ converges to
$\zero$.

Let $[B,\varepsilon]^s(\zero)$ be a neighbourhood of $\zero$.
Since $\{U_{n,n}:n\in \naturals\}\in \Gamma_{\mathfrak B^s}$,
there are $m\in\naturals$ and $\delta_n >0$, $n > m$, such that
$1/m < \varepsilon$ and for each $n> m$, $B^{2\delta_n} \subset
U_{n,n}$. Therefore, for all $n>m$, $f_{B_n,n,n}(B^{\delta_n})
\subset (-1/n,1/n) \subset (-\varepsilon,,\varepsilon)$, i.e.
$f_{B_n,n,n} \in [B,\varepsilon]^s(\zero)$. $\eproof$

\begin{theorem} \label{s1-sub} If $(X,d)$ is a metric space and $\mathfrak B$ a
bornology with closed base on $X$, then the following statements are
equivalent:
\begin{itemize}
\item[$(1)$] $X$ satisfies $\sone(\mathcal O_{\mathfrak
B^s},\Gamma_{\mathfrak B^s})$;
\item[$(2)$] Each $\mathfrak B^s$-cover $\mathcal U$ of $X$ contains a countable set
$\{U_n:n\in \naturals\}$ which is a $\gamma_{\mathfrak B^s}$-cover
of $X$.
\end{itemize}
\end{theorem}
$\proof$ Obviously $(1) \Rightarrow (2)$ and thus we prove $(2)
\Rightarrow (1)$.

$\proof$ $(2)\Rightarrow (1)$: [Observe first that each space
satisfying (2) is $\mathfrak B^s$-Lindel\"of.] Let $\{\mathcal U_n :
n \in \naturals\}$ be a sequence of $\mathfrak B^s$-covers of $X$.
Construct a new sequence $(\mathcal V_n:n\in \naturals)$ of
$\mathfrak B^s$-covers of $X$ in the following way:

(i) $\mathcal V_1 = \mathcal U_1$;

(ii) $\mathcal V_{n+1}$ is a refinement of $\mathcal V_n$ and
$\mathcal U_n$.

As the bornology $\mathfrak B$ has closed base, by Fact 1, for each
$B\in \mathfrak B$ the closure of $B$ is also in $\mathfrak B$. On
the other hand,  $X\notin \mathfrak B$, so that for every $B\in
\mathfrak B$ there is a point $x_B \in X \setminus \overline{B}$. It
follows there is $\delta>0$ such that $x_B\notin B^\delta$, i.e.
$B^\delta \subset X\setminus \{x_B\}$. Therefore, $\{X\setminus
\{x\}:x\in X\}$ is a $\mathfrak B^s$-cover of $X$. Applying
assumption (2) to $\{X\setminus \{x\}:x\in X\}$ we find a sequence
$(x_n:n\in\naturals)$ of points of $X$ such that $\{X\setminus
\{x_n\}:n\in \naturals\}$ is still a $\mathfrak B^s$-cover of $X$. For
each $n \in \naturals$ denote
\[
\mathcal W_n = \{V \setminus \{x_n\} : V \in \mathcal V_n\} \mbox{
and } \mathcal W = \cup\{ \mathcal W_n : n \in \naturals\}.
\]
We claim that $\mathcal W$ is a $\mathfrak B^s$-cover of $X$.

Let $B\in \mathfrak B$. There are $\delta > 0$ and $k\in\naturals$
such that $B^{\delta} \subset X\setminus \{x_k\}$. Also, since
$\mathcal V_k$ is a $\mathfrak B^s$-cover, there are $\mu >0$ and
$V\in \mathcal V_k$ such that $B^{\mu} \subset V$. Then for
$\varepsilon = \min\{\delta,\mu\}$ we have $B^{\varepsilon} \subset
V\setminus \{x_k\} \in \mathcal W_k \subset \mathcal W$.

By $(2)$, there exists a sequence $(W_m: m \in \naturals)$ in
$\mathcal W$ so that $\{W_m:m\in\naturals\} \in \Gamma_{\mathfrak
B^s}$. For every $m \in \naturals$ there exist $n_m \in \naturals$
and a set $V_{n_m} \in \mathcal V_{n_m}$ with $W_m \subset
V_{n_m}$. Let for each $m \in \naturals$, $F_m = \{x_1,...,
x_m\}$; clearly, $F_m\in \mathfrak B$.

Choose $W_{m_1} \in \{ W_m : m \in \naturals\}$ such that
$F_1^{\delta} \subset W_{m_1}$ for some $\delta >0$; note that
$m_1
> 1$. For each $p > 1$ pick $W_{m_p} \in \{ W_m : m \in \naturals\}$ so that
$F_p^{\mu} \subset W_{m_p}$ for some $\mu >0$, $m_p > m_{p-1}$ and
$m_p > p$. Such a choice is possible for each $p \in \naturals$
because the set $\{ W_m : m \in \naturals\}$ is a $\gamma_{\mathfrak
B^s}$-cover of $X$. Since $\mathcal V_{n+1}$ is a refinement of both
$\mathcal V_n$ and $\mathcal U_n$, for each $i \in \naturals$ with
$n_{m_p}< i\le n_{m_{p+1}}$ pick $U_i \in \mathcal U_i$ such that
$U_i \subset U_j$ for $i>j$. Put $n_{m_0}=0$ and for each $n \in
\naturals$ define (i)
$O_n = U_{n_{m_p}}$ if $n = n_{m_p}$, and (ii) $O_n = U_n$ if $n_{m_p}< n<
n_{m_{p+1}}$.

We claim that the sequence $(O_n : n \in \naturals)$ is a selector
for $(\mathcal U_n : n \in \naturals)$ witnessing that $X$
satisfies $\sone (\mathcal O_{\mathfrak B^s}, \Gamma_{\mathfrak
B^s})$. Indeed, let $B \in \mathfrak B$. Then there exist $m_0 \in
\naturals$ and a sequence $(\delta_n:n > m_0)$ such that for every
$m> m_0$, $B^{\delta_n} \subset U_m$. By construction of the sets
$O_n$, for all $n> n_{m_0}$ we have $B^{\delta_n} \subset O_n$.
$\eproof$

\medskip
In \cite{agataholapeppe} it was proved:

\begin{theorem} {\rm (\cite{agataholapeppe})} If $(X,d)$ is a metric space
and $\mathfrak B$ a bornology with closed base on $X$, then the
following are equivalent:
\begin{itemize}
\item[$(1)$] $({\sf C}(X),\tau_{\mathfrak B}^s)$ is a Fr\'echet-Urysohn space;
\item[$(2)$] Each $\mathfrak B^s$-cover $\mathcal U$ of $X$ contains a countable set
$\{U_n:n\in \naturals\}$ which is a $\gamma_{\mathfrak B^s}$-cover
of $X$.
\end{itemize}
\end{theorem}

From this theorem together with Theorem \ref{SFU} and Theorem \ref{s1-sub}
one obtains the following corollary. (Similar well-known results were obtained independently
in \cite{gerlitsnagy} and \cite{pytkeev} for the space ${\sf
C}_p(X)$ (see \cite{arhbook}). Compare also with \cite{kcov2},
\cite{chinese}, \cite{mccoy} in connection with the space ${\sf
C}_k(X)$.)

\begin{corollary} \label{FU-SFU} If $(X,d)$ is a metric space
and $\mathfrak B$ a bornology with closed base on $X$, then the
following assertions are equivalent:
\begin{itemize}
\item[$(1)$] $({\sf C}(X),\tau_{\mathfrak B}^s)$ is a Fr\'echet-Urysohn space;
\item[$(2)$] $({\sf C}(X),\tau_{\mathfrak B}^s)$ is a strictly Fr\'echet-Urysohn space;
\item[$(3)$] Each $\mathfrak B^s$-cover $\mathcal U$ of $X$ contains a countable set
$\{U_n:n\in \naturals\}$ which is a $\gamma_{\mathfrak B^s}$-cover
of $X$;
\item[$(4)$] $X$ satisfies $\sone(\mathcal O_{\mathfrak
B^s},\Gamma_{\mathfrak B^s})$.
\end{itemize}
\end{corollary}

We close this section considering a property similar to the strict Fr\'echet-Urysohn property.

A space $X$ is called a \emph{selectively strictly $A$-space} (shortly ${\sf SSA}$) \cite{vlada}
if for each sequence $(A_n:n\in\naturals)$ of subsets
of $X$ and each point $x\in X$ such that $x\in \overline
{A_n}\setminus A_n$ for each $n\in \naturals$, there is a sequence
$(T_n:n\in \naturals)$ such that for each $n$ \, $T_n\subset A_n$, and
$x\in\overline{\bigcup_ {n\in\naturals}T_n}
\setminus\bigcup_{n\in\naturals}\overline{T_n}$.

\begin{theorem}\label{selAspace}
Let $(X,d)$ be a metric space and $\mathfrak B$ a bornology with closed
base on $X$. Then the following are equivalent:
\begin{itemize}
\item[$(1)$] $({\sf C}(X),\tau_{\mathfrak B}^s)$ is ${\sf SAA}$;
\item[$(2)$] for each sequence $(\mathcal U_n:n\in\naturals)$ of $\mathfrak B
^s$-covers of $X$ there is a sequence $(\mathcal V_n:n\in\naturals)$ such that
$\mathcal V_n \subset \mathcal U_n$ for each $n$, no $\mathcal V_n$
is a $\mathfrak B^s$-cover of $X$, and $\bigcup_{n\in\naturals}\mathcal
V_n$ is a $\mathfrak B^s$-cover of $X$.
\end{itemize}
\end{theorem}
$\proof$ $(1) \Rightarrow (2)$: Let $(\mathcal U_n:n\in\naturals)$
be a sequence of $\mathfrak B^s$-covers of $X$. For every
$n\in \naturals$ and every $B\in \mathfrak B$ there are $U_B \in
\mathcal U_n$ and $\delta >0$ with $B^{2\delta} \subset U_B$.

For $n\in\naturals$, $B\in\mathfrak B$ put
\[
\mathcal{U}_{n,B}:=\{U\in{\mathcal U} _n: B^{2\delta}\subset U\}.
\]
For each $U\in\mathcal{U}_{n,B}$ pick a continuous
function $f_{n,B,U}:X\to [0,1]$ satisfying
$f_{n,B,U}(B^{\delta})= \{0\}$ and $f_{n,B,U}(X\setminus U)=\{1\}$ and denote
\[
A_n=\{f_{n,B,U}:B\in  \mathfrak B, U\in \mathcal{U}_{n,B}\} \, \, (n\in \naturals).
\]
The function $\zero$ belongs to the $\tau_{\mathfrak B}^s$-closure
of $A_n$ as it is easy to verify. On the other hand, $\zero \notin
A_n$ for each $n\in\naturals$. Otherwise $\zero = f_{m,B,U} \in
A_m$ for some $m$, thus $f_{m,B,U}(X\setminus U)= \{0\}$, a
contradiction.

By (1), there is a sequence $(T_n:n\in \naturals)$ such that for each
$n$ \, $T_n\subset A_n$ and $\zero \in\overline{\bigcup_ {n\in\naturals}T_n}
\setminus \bigcup_{n\in\naturals}\overline{T_n}$. Denote by $\mathcal V_n$, $n\in \naturals$,
the set of corresponding sets $U$ for each $f_{n,B,U}\in T_n$. No $\mathcal V_n$ is a $\mathfrak
B^s$-cover of $X$, because otherwise $\zero \in \overline{T_n}$.
It remains to prove $\bigcup_{n\in\naturals}\mathcal V_n \in \mathcal O_{\mathfrak B^s}$.
Let $B\in \mathfrak B$. From $\zero \in \overline{\bigcup_ {n\in\naturals}T_n}$ it follows the existence
of $m\in\naturals$ and $f_{m,B,U} \in T_m$ such that there is $\delta >0$ with $f_{m,B,U}(x) = 0$ for each
$x\in B^{\delta}$. This means $B^{2\delta} \subset U\in \mathcal V_m \subset \bigcup_{n\in\naturals}\mathcal V_n$.\\

$(2) \Rightarrow (1)$: Let $(A_n:n\in\naturals)$ be a sequence of subsets of
$({\sf C}(X),\tau_{\mathfrak B}^s)$ such that $\zero\in
\overline{A_n} \setminus A_n$, $n\in
\naturals$. We proceed as in the proofs of $(2)\Rightarrow (1)$ of Theorems \ref{strongfan} and \ref{SFU}.
For each $B\in \mathfrak B$
and each $m\in \naturals$ the neighborhood $[B,1/m]^s(\zero)$ of
$\zero$ intersects each $A_n$. Therefore, for each $n\in \naturals$
there is a function $f_{B,n,m}\in A_n$ such that there is $\delta
> 0$ with $|f_{B,n,m}(x)|< 1/m$ for each $x\in B^{\delta}$. For
each $n,m \in \naturals$ define
\[
{\mathcal U}_{n,m}=\{f^{\gets}(-1/m,1/m):m\in \naturals, f\in A_n\}.
\]
All ${\mathcal U}_{n,m}$ are $\mathfrak B^s$-covers of $X$.

\smallskip
Case 1: $X\in \mathcal U_{n,n}$ for infinitely many $n$.

Then there exist an increasing sequence $n_1<n_2<\dots<n_k<\dots$ in
$\naturals$ and $f_{n_k}\in A_{n_k}$, $k\in\naturals$, such that
$f_{n_k}^{\gets}(-1/n_k, 1/n_k)=X$. Put $T_{n_k}=\{f_{n_k}\}$,
$k\in\naturals$, and $T_{n}=\emptyset$ for $n\neq n_k$,
$k\in\naturals$. It is easy to see that $\zero \in \overline{\bigcup_{n\in
\naturals}T_n}\setminus \bigcup_{n\in \naturals}\overline{T_n}$ (because the sequence
$(f_{n_k}:k\in \naturals)$ actually $\tau_{\mathfrak B}^s$-converges to $\zero$).

\smallskip
Case 2: $X\in\mathcal U_{n,n}$ for finitely many $n$.

Suppose that $X\notin \mathcal U_{n,n}$ for each $n$. Choose a sequence
$(\mathcal V_{n,n}:n\in\naturals)$ witnessing (2).
For each $V\in \mathcal V_{n,n}$ pick a function $f_V\in A_n$ with $V= f_V^{\gets}
(-1/n, 1/n)$ and put $T_n =\{f_V:V\in \mathcal V_{n,n}\}$. Let us show
that $(T_n:n\in\naturals)$ is as required in (1).

First, $\zero \notin \bigcup_{n\in \naturals}\overline{T_n}$. Otherwise,
$\zero \in \overline{T_m}$ for some $m$ would imply that
$\mathcal V_{m,m}$ is a $\mathfrak B^s$-cover of $X$.

Second, we prove $\zero \in \overline{\bigcup_{n\in\naturals}T_n}$.
Let $[B,\varepsilon]^s(\zero)$, $B\in\mathfrak B$, $\varepsilon >0$,  be a neighborhood
of $\zero$. Since $\bigcup_{n\in\naturals} \mathcal V_{n,n}$ is a
$\mathfrak B^s$-cover of $X$ there are $\delta > 0$ and a natural number $m$ such
that $1/m < \varepsilon$
and for some $V\in \mathcal V_{m,m}$ we have $B^{\delta} \subset V =
f_V^{\gets}(-1/m, 1/m)$. Hence $f_V\in T_m$ and since obviously $f_V\in [B,\varepsilon]^s(\zero)$,
we conclude $\zero \in \overline{\bigcup_{n\in\naturals}T_n}$.
$\eproof$


\section{Covering properties}

In this short section we indicate how the idea of strong uniform convergence
on a bornology may be further applied to selection principles theory.

According to \cite{taiwanese}, we say that a space $X$ with a bornology $\mathfrak B$
satisfies the selection principle $\alpha_i
(\Gamma_{\mathfrak B^s},\Gamma_{\mathfrak B^s})$, $i = 2,3,4$, if
for each sequence $(\mathcal U_n:n\in\naturals)$ of $\gamma_{\mathfrak B^s}$-covers
of $X$ there is a $\gamma_{\mathfrak B^s}$-cover $\mathcal V$ of $X$ such that:

$\alpha_2(\Gamma_{\mathfrak B^s},\Gamma_{\mathfrak B^s})$: for each
$n\in\naturals$, $\mathcal U_n \cap \mathcal V$ is infinite;

$\alpha_3(\Gamma_{\mathfrak B^s},\Gamma_{\mathfrak B^s})$: for infinitely many
$n\in\naturals$, $\mathcal U_n \cap \mathcal V$ is infinite;

$\alpha_4(\Gamma_{\mathfrak B^s},\Gamma_{\mathfrak B^s})$: for infinitely many
$n\in\naturals$, $\mathcal U_n \cap \mathcal V \neq \emptyset$.

\begin{theorem} \label{alpha234-gamma-k} For a metric space $(X,d)$
and a bornology $\mathfrak B$ on $X$ the following are equivalent:
\begin{itemize}
\item[$(1)$] $X$ satisfies $\alpha_2(\Gamma_{\mathfrak B^s},
\Gamma_{\mathfrak B^s})$;
\item[$(2)$] $X$ satisfies
$\alpha_3(\Gamma_{\mathfrak B^s}, \Gamma_{\mathfrak B^s})$;
\item[$(3)$]
$X$ satisfies $\alpha_4(\Gamma_{\mathfrak B^s}, \Gamma_{\mathfrak
B^s})$;
\item[$(4)$] $X$ satisfies $\sone(\Gamma_{\mathfrak B^s},
\Gamma_{\mathfrak B^s})$.
\end{itemize}
\end{theorem}
$\proof$ Clearly $(1) \Rightarrow (2) \Rightarrow (3)$.

\smallskip
$(3) \Rightarrow (4)$: Let $(\mathcal U_n:n\in\naturals)$ be a
sequence of $\gamma_{\mathfrak B^s}$-covers of $X$. Suppose that
$\mathcal U_n = \{U_{n,m}:m\in \naturals\}$, $n\in\naturals$. For
all $n,m\in\naturals$ define $V_{n,m} = U_{1,m}\cap U_{2,m} \cap
\cdots \cap U_{n,m}$. Then for each $n$ the set $\mathcal V_n =
\{V_{n,m}:m\in \naturals\}$ is a $\gamma_{\mathfrak B^s}$-cover of
$X$. (Indeed, fix $n$ and take a $B\in\mathfrak B$.  Since each
$\mathcal U_i$, $i\le n$, is a $\gamma_{\mathfrak B^s}$-cover of
$X$, there are $m_i$ and sequences $(\delta_m^{(i)}:m\ge m_i)$,
$i\le n$, of positive real numbers, such that
$B^{\delta_m^{(i)}}\subset U_{i,m}$ for every $m\ge m_i$.  The
sequence $(\delta_m^{(i)}:i \le n, m\ge m_0)$ and $m_0 =
\max\{m_i:i\le n\}$ witness that $\mathcal V_n$ is a
$\gamma_{\mathfrak B^s}$-cover.) By (3) (and the fact that an
infinite subset of a $\gamma_{\mathfrak B^s}$-cover is also such a
cover) there is an increasing sequence $n_0 = 1 \le n_1 < n_2 <
\cdots$ in $\naturals$ and a $\gamma_{\mathfrak B^s}$-cover
$\mathcal V = \{V_{n_i,m_i}:i\in\naturals\}$ such that for each
$i\in \naturals$, $V_{n_i,m_i}\in \mathcal V_{n_i}$. For each $i
\ge 0$, each $j$ with $n_i < j \le n_{i+1}$  and each
$V_{n_{i+1},m_{i+1}} = U_{1,m_{i+1}} \cap \cdots \cap
U_{n_{i+1},m_{i+1}}$ let $U_j$ be the set $U_{j,m_{i+1}}$. The
sequence $(U_n:n\in\naturals)$ testifies that $X$ satisfies
$\sone(\Gamma_{\mathfrak B^s}, \Gamma_{\mathfrak B^s})$.\\

$(4) \Rightarrow (1)$: Let $(\mathcal U_n: n\in\naturals)$ be a
sequence of $\gamma_{\mathfrak B^s}$-covers of $X$. Suppose as
above $\mathcal U_n = \{U_{n,m}:m\in \naturals\}$,
$n\in\naturals$. As we mentioned before every infinite subset of
a $\gamma_{\mathfrak B^s}$-cover is itself a $\gamma_{\mathfrak
B^s}$-cover. Thus from each $\mathcal U_n$ we can built countably
many disjoint $\gamma_{\mathfrak B^s}$-covers $\mathcal U_n^k$,
$k\in \naturals$. Apply now (4) to the sequence $(\mathcal
U_n^k:n,k\in \naturals)$. One obtains a $\gamma_{\mathfrak
B^s}$-cover $\{U_n^k:U_n^k\in \mathcal U_n^k; n,k \in \naturals\}$
which contains infinitely many elements from each $\mathcal U_n$.
$\eproof$

\medskip
Observe that similar results are true for the pairs  $(\Gamma_s,
\Gamma_s)$, $(\mathcal O_{\mathfrak B^s}, \Gamma_s)$ and $(\mathcal
O_{\mathfrak B^s}, \Gamma_{\mathfrak B^s})$.

\vspace{5mm} \footnotesize{

}

\medskip\normalsize
\begin{center}
{\bf Addresses}
\end{center}

\begin{center}
\begin{tabular}{ll}
Agata Caserta                              & Giuseppe Di Maio \\
Dipartimento di Matematica                 & Dipartimento di Matematica \\
Seconda Universit$\grave{a}$ di Napoli     & Seconda Universit$\grave{a}$ di Napoli \\
Via Vivaldi 43                             & Via Vivaldi 43 \\
81100 Caserta                              & 81100 Caserta \\
Italy                                      & Italy \\
{\sf agata.caserta@unina2.it}              & {\sf giuseppe.dimaio@unina2.it} \\
                                           & \\

Ljubi\v sa D.R. Ko\v cinac\\
Faculty of Sciences and Mathematics \\
University of Ni\v s\\
Vi\v segradska 33\\
18000 Ni\v s\\
Serbia\\
{\sf lkocinac@gmail.com} \\
\end{tabular}
\end{center}

\end{document}